\documentclass[11pt,
]{article}

\usepackage[tbtags]{amsmath}
\usepackage{amsfonts,amssymb,mathrsfs,amscd,amstext,
latexsym,amsbsy,amsthm}

\usepackage{mparhack}
\usepackage{marginnote}
\usepackage{xcolor}

\usepackage[bbgreekl]{mathbbol}

\input{axioms.sty}

\vyk{
\nek{\toti}{%
\fontfamily{ptm}\selectfont}
}

\usepackage{doi}

\begin{document}

%\widowpenalty-10

\title
{A model in which the Separation principle holds 
for a given effective  
projective Sigma-class\thanks
{This research was partially supported by 
Russian Foundation for Basic Research RFBR grant number 
20-01-00670.}
}

\author 
{Vladimir Kanovei\thanks{
IITP RAS,
 \ {\tt kanovei@iitp.ru} .
}  
\and 
Vassily Lyubetsky\thanks{IITP RAS,
\ {\tt lyubetsk@iitp.ru}. 
%Partial support of grant RFBR 18-29-13037 acknowledged.
} 
}

\date 
{\today}

\maketitle

%$\,$\hfill ÓÄÊ 510.225 è 510.223

\begin{abstract}
In this paper, we prove the following: 
If $n\ge3$, there is a generic extension 
of $\rL$ -- the constructible universe --  
in which it is 
true that the Separation principle holds for both 
effective (lightface) classes 
%\highlighting
{$\is1n$} 
 and 
%\highlighting
{$\ip1n$ }
of sets of integers.  
The result was announced long ago by Leo Harrington 
with a sketch of the proof for $n=3$; 
its full proof has never been presented.
Our methods are based on a countable product of 
almost-disjoint forcing notions independent in the sense 
of Jensen--Solovay.
\end{abstract}

\parf{Introduction}
\las{int}

The separation problem was introduced in 
descriptive set theory by Luzin~\cite{lbook}. 
In modern terms, the separation principle -- or 
simply 
%\hl
{\sep}, for a given projective (boldface) 
class $\fs1n$ or $\fp1n$ -- is the assertion that\vom\vim 
\bde
\item[\bsep\ 
%\hl
{\rm for} 
%MDPI: please check this part format, is it list format?
%%%% IT'S `description' FORMAT
%%%% `for' AND `or' ARE MADE NON-BOLD
 $\fs1n$ {\rm or} $\fp1n$:] 
any pair of disjoint 
$\fs1n$, resp,\ $\fp1n$ sets $X,Y$ of reals 
can be separated by a $\fd1n$ set.\vom\vim 
\ede
Accordingly, 
\rit{the classical separation problem} 
%Is the italics necessary? 
%%%% YES IT IS
is a question of whether 
\bsep\ holds for this or another projective class $\fs1n$ 
or $\fp1n$.     
Luzin and then Novikov \cite{nov1951} underlined the 
importance and difficulty of this problem. 
(See \cite{mDST,Kdst,kl49} for details and further references.)

Luzin \cite{lus:ea,lbook} and Novikov~\cite{nov1931} 
proved that \bsep\ holds  
for $\fs11$ but fails for the dual class $\fp11$.
Somewhat later, it was established 
by Novikov~\cite{nov1935} that the picture changes 
at the next projective level: 
\bsep\ holds for $\fp12$ but fails for $\fs12$. 

As for the higher levels of projective hierarchy, 
all attempts made in classical descriptive set theory 
to solve the separation problem above the 
second level did not work, 
until some additional set theoretic axioms were 
introduced---in particular, those by Novikov~\cite{nov1951} 
and Addison~\cite{add1,add2}. 
G\"odel's \rit{axiom of constructibility} 
$\rV=\rL$ implies that, for any $n\ge 3$, 
\bsep\ holds for $\fp1n$ but fails for 
$\fs1n$---pretty similar to second level.

In such a case, it is customary in modern set theory to 
look for models in which the separation problem is solved 
differently than under $\rV=\rL$ for at least 
some projective classes $\fs1n$ and $\fp1n$, $n\ge3$. 
This goal is split into two different tasks: 
\ben
\Renu
\itlb{tas1}% 
Prove the \rit{independence} of the \dd\Pi side 
\bsep---that is, given $n\ge3$, 
find models in which \bsep\ \rit{fails} for the 
class $\fp1n$;\vim

\itlb{tas2}%
Prove the \rit{consistency} of the \dd\Sigma side 
\bsep---that is, given $n\ge3$, 
find models in which \bsep\ \rit{holds} 
for the class $\fs1n$.\vom\vim
\een
As for \rit{models}, we focus here only on 
\rit{generic extensions of 
the constructible universe $\rL$}. 
Other set theoretic models, e.g.,
those based on strong determinacy 
or large cardinal hypotheses, 
%as in Footnote~\ref{sno1}, 
are not considered in this paper. 
({%\hl%
We may only note in brackets that, 
by Addison and Moschovakis   
\cite{addmos}, 
and Martin~\cite{martAD}, 
the \rit{axiom of projective determinacy} $\PD$  
implies that, for any $m\ge 1$, 
the separation problem is solved affirmatively  
for $\fs1{2m+1}$ and $\fp1{2m+2}$ and   
negatively for $\fp1{2m+1}$ and $\fs1{2m+2}$ 
--- similar to what happens at the first and second level 
corresponding to $n=0$ in this scheme. 
See also Steel~\cite{steel_detsep,steel_core}, and
Hauser and Schindler \cite{hashi} for some other 
relevant results.}).

%\hl
{Problems} 
(I) and (II) have been well-known 
since the early years of forcing, e.g., see problem
P3030, and especially P3029 
(= (II) for $n=3$) 
in a survey \cite{mathiasS} by Mathias. 

Two solutions for part (I) are known so far. 
Harrington's two-page handwritten note  
([Addendum A1]~\cite{h74}) 
contains a sketch of a model, 
without going into details, defined by the technique of 
almost-disjoint forcing of Solovay and Jensen \cite{jsad}, 
in which indeed \sep\ fails for  
both $\fs1n$ and $\fp1n$ for a given $n$. 
This research was cited in 
%MDPI: revised this format, please check it.
%%%% WE HAVE IT SLIGHTLY CORRECTED
Moschovakis~\cite{mDST} (Theorem 5B.3),    
and Mathias~\cite{mathiasS} (Remark P3110 on page 166), 
but 
has never been published or otherwise detailed in any way. 
Some other models, with the same property 
of failure of \sep\ for different projective classes,  
were recently defined and studied in \cite{kl28,kl49}. 

As for (II), the problem as it stands is open 
so far, and no conclusive achievement, such as a model 
(a generic extension of $\rL$) 
in which \bsep\ holds for $\fs1n$ for some $n\ge3$, 
is known. 
Yet, the following modification turns out to be 
easier to work with. 
The \rit{effective} or \rit{lightface} \sep, for a 
given lightface class $\is1n$ or $\ip1n$ 
(we give \cite{mDST} as a reference on the lightface 
projective hierarchy), 
is the assertion that\vom\vim 

\bde
\item[\lsep\ {\rm for} $\is1n$ {\rm or} $\ip1n$:] 
any pair of disjoint 
$\is1n$, respectively,\ $\ip1n$ sets $X,Y$ 
%of reals 
can be separated by a $\id1n$ set---here,
unlike the \bsep\ case, the sets $X,Y$ can be either 
sets of reals or sets of integers.\vom\vim 
\ede
Accordingly, the \rit{effective} or \rit{lightface} 
separation problem is a question of whether   
\lsep\ holds for this or another class of the form 
$\is1n$ or $\ip1n$, with specific versions for 
sets of reals and sets of integers.  
Addison~\cite{add1,add2} demonstrated that, 
similar to the above, \lsep\ 
holds for $\is11$ and $\ip12$; 
fails for $\ip11$ and $\is12$; 
and under the axiom of constructibility $\rV=\rL$, it  
holds for $\ip1n$ and fails for $\is1n$ for all 
$n\ge3$---both in the ``real'' and the ``integer'' versions. 
(See also \cite{mDST}.) 

In this context, Harrington announced 
in \cite{h74} that there is a model in which 
\lsep\ holds both for the class $\is13$ 
\rit{for sets of integers}, 
and for the class $\ip13$ 
\rit{for sets of integers}. 
A two-page handwritten sketch of a construction of 
such a model is given in ([Addendum A3]~\cite{h74})
without much elaboration of arguments. 
{\ubf The goal of this paper} is to prove the next 
theorem, 
which generalizes the cited Harrington result and 
thereby is a definite advance in the direction of 
(II) in the context of \lsep\ for 
sets of integers. 
This is the main result of this paper.

%\begin{Theorem}
%
\bte
\lam{mt}
Let\/ $\nn\ge2$. 
There is a generic extension of\/ $\rL$ in which\/ 
\ben
\renu
\jtlb{mt1}%
\lsep\ holds for\/ $\is1{\nn+1}$ sets of integers, 
so that 
any pair of disjoint\/ $\is1{\nn+1}$ sets\/ $X,Y\sq\om$  
can be separated by a\/ $\id1{\nn+1}$ set$;$

\jtlb{mt2}%
\lsep\ also holds for\/ $\ip1{\nn+1}$ sets of integers, 
so that 
any pair of disjoint\/ $\ip1{\nn+1}$ sets\/ $X,Y\sq\om$  
can be separated by a\/ $\id1{\nn+1}$ set.
\een
\ete

%[Addendum A3]

Our proof of this theorem 
will follow a scheme that includes both 
some arguments outlined by Harrington in {\cite{h74}}, 
Addendum A3
(mainly related to the most elementary case $\nn=2$) 
and some arguments absent in 
{\cite{h74}}, 
in particular, those related to the generalization to 
the case $\nn\ge3$. 
(We may note here that {\cite{h74}} is neither a 
beta-version of a paper, nor a preprint of any sort, 
but rather handwritten notes to a talk in which omissions 
of even major details can be expected.) 
All this will require both a fairly sophisticated 
construction of the model itself and a fairly complex 
derivation of its required properties by rather 
new methods for modern set theoretic research. 
Thus, we are going to proceed with filling-in all 
necessary details left aside in {\cite{h74}}. 
We hope that the detailed acquaintance with the 
set theoretic methods first introduced by Harrington  
will serve to the benefits of the reader~envisaged.

To prove Theorem \ref{mt}, we make use of a 
generic extension of $\rL$ defined in our 
earlier paper \cite{kl57} 
(and before that in \cite{h74}---modulo some key details absent in \cite{h74}) 
in order to prove the consistency of the equality 
$\pws\om\cap\rL=\pws\om\cap\id1{\nn+1}$ 
for a given $\nn\ge2$. 
(The equality 
claims that the constructible reals are the same as 
$\id1{\nn+1}$ reals. 
Its consistency was a major problem posed by  
Harvey Friedman \cite{102}.) 
We present the construction of this generic model in 
all necessary detail. 

This includes a version of almost-disjoint forcing 
considered in Section~\ref{prel}, the cone-homogeneity 
lemma in Section~\ref{fhom}, the systems and 
product forcing construction in Section~\ref{*s41}, 
and a Jensen--Solovay-style construction of the actual 
countable support forcing product $\dQ$ in Section~\ref{*71}.  
Theorem~\ref{*7euv} and Definition~\ref{*fixu} 
in Section~\ref{*71} present the construction of $\dQ$ 
in $\rL$ via the union of a \dd\omd long increasing 
sequence of systems $\xU_\xi$, $\xi<\omb$, 
which satisfies  
suitable completeness and definability requirements 
(that depend on the choice of the value of an integer 
$\nn$ as in Theorem~\ref{mt}),  
and also follows the Jensen--Solovay  
idea of Cohen-generic 
extensions at each step $\xi<\omb$ of the inductive 
construction of $\xU_\xi$. 

Then, we consider 
corresponding \dd\dQ generic extensions $\rL[G]$
in Section~\ref{*Pex}, 
and their subextensions involved in the proof 
of Theorem~\ref{mt} in Section~\ref{175} 
(Theorem~\ref{33}). 
Two key lemmas are established in Section~\ref{2key}, 
and the proof of theorems \ref{33} and \ref{mt} 
is finalized in 
Section~\ref{175*} (the $\is1{\nn+1}$ case) and in 
Section~\ref{176} (the $\ip1{\nn+1}$ case). 

The final section on conclusions and discussion 
completes the paper.

%On the other hand, we'll skip some proofs, giving 
%reference to the corresponding results fully proved 
%in \cite{kl57}. 

%%%%%%%%%%%%%%%%%%%%%%%%%%%%%%%%%%%%%%%%%%
\parf{Almost-disjoint forcing}
\las{prel}

Almost-disjoint forcing 
was invented by Jensen and Solovay~\cite{jsad}. 
Here, we make use of a \dd\omi version of this tool 
considered in ([Section 5]~\cite{jsad}). 
The version we utilize here exactly corresponds to 
the case $\bom=\omil$ developed 
in our earlier paper~\cite{kl57} 
(and to less extent in \cite{kl58}). 
This will allow us to skip some proofs below. 
Assume the following: %please confirm if intended meaning is retained.
\bit
\item\msur
$\bom=\omil$, the first uncountable ordinal in $\rL$.

%\vitem\msur
%$\bomp=\om^\rL_2$, the 2nd uncountable ordinal in $\rL$;

\vitem\msur
$\Fun=({\bom}^{\bom})\cap\rL$ = all \dd\bom sequences 
of ordinals $<\bom$ in $\rL$. 

\vitem\msur
$\Seq=({\bom}^{<\bom}\bez\ans\La)\cap\rL$ 
is the set of  
all sequences $s\in\rL$ 
of ordinals $<\bom$, of length $0<\lh s<\bom$.
\eit

By definition, the sets 
$\Fun$, $\Seq$ belong to $\rL$ and 
$\card{(\Seq)}=\bom=\omil$ whereas 
$\card{(\Fun)}=\om^\rL_2$ in $\rL$.
Note that $\La$, the~empty sequence, 
does not belong to $\Seq$.
\bit
\item 
A set $X\sq\Fun$ is \rit{dense} iff for 
any $s\in\Seq$ there is 
\index{set!dense}%
$f\in X$ such that $s\su f$.

\vitem
If $S\sq\Seq\yt f\in \Fun$; then, let 
$\ssl{S}f= %\tsup
\ens{\xi<\bom}{f\res\xi\in S}$.
\index{zzsf@$\ssl{S}f$}%

\vitem
If $\ssl{S}f$  
\imar{ssl Sf}%
\index{zzsf@$\ssl S f$}%
is unbounded in $\bom$, then say that
\index{cover}%
\rit{$S$ covers\/ $f$}; otherwise,  
\rit{$S$ does not cover\/ $f$}.
\eit

The general goal of the almost-disjoint forcing 
is the following:  
given a set $u\sq\Fun$ in the ground set universe $\rL$, 
find a generic set $S\sq\Seq$ such that 
the equivalence
``$f\in u\leqv S\:\text{ does not cover }f$'' 
holds for each $f\in \Fun$ in the ground universe. 
\bdf
%[in $\rL$]
\lam{*3.2a}
\imar{usl ps p pf p}%
$\Qa$ is the set of all pairs\/ 
$p=\usl{\ps p}{\pf p}\in\rL$
 of {\em finite} 
\index{zzQ*@$\Qa$}%
\index{zzsp@$\ps p$}%
\index{zzfp@$\pf p$}%
sets\/ $\pf p\sq\Fun$, $\ps p\sq\Seq$. 
Note that\/ $\Qa\in\rL$. 
Elements of\/ $\Qa$ are called 
(forcing) \rit{conditions}. 
\index{condition}%

If\/ $p\in\Qa$, then put\/
\imar{pfv p}%
$
\pfv p=\ens{f\res\xi}{f\in\pf p\land 1\le\xi<\bom}\,,
$; 
this is a tree in $\Seq$. 

Let\/ $p,q\in\Qa$.
Define\/ $q\leq p$ 
(that is, $q$ is {\em stronger} as a forcing condition) 
iff\/ 
\index{condition!order}%
\index{condition!stronger}%
$\ps p\sq \ps q\yd \pf p\sq \pf q$, 
and~the difference\/ $\ps q\bez \ps p$ 
does not intersect $\pfv p$, \ie,
$\ps q\cap\pfv p=\ps p\cap \pfv p$.
Clearly, we have\/ $q\leq p$ iff\/  
$\ps p\sq \ps q\yd \pf p\sq \pf q$, 
and\/  
$\ps q\cap {\pfv p}=\ps p\cap {\pfv p}$.  

If $u\sq\Fun$, then put\/ 
$\zq u=\ens{p\in \Qa}{\pf p\sq u}$.
\index{zzQu@$\zq u$}%
\imar{zq u}
\edf

\ble
[Lemma 1 in \cite{kl57}]
\lam{*comp}
Conditions\/ $p,q\in\Qa$ are compatible in\/ $\Qa$ iff\/ 
$(1)$ $\ps q\bez \ps p$ does not intersect\/ $\pfv p$, and\/ 
$(2)$ $\ps p\bez \ps q$ does not intersect\/ $\pfv q$.
%\hl{\qed} 
%MDPI: suggest remove it.
%%%% REMOVED
\ele

%\ele

%\bcor
%\lam{*cor}
Thus, any conditions\/ $p,q\in \zq u$ are compatible 
in\/ $\zq u$ iff\/ $p,q$ are compatible in\/ $\Qa$ iff 
the condition\/ 
$p\land q=\usl{\ps p\cup \ps q}{\pf p\cup \pf q}\in \zq u$ 
\index{zzpandq@$p\land q$}%
satisfies\/ $p\land q\leq p$ and\/ $p\land q\leq q$. 
%\qed
%\ecor

%

\parf{The almost-disjoint forcing notions 
are homogeneous}
\las{fhom}

We are going to show that forcing notions of the 
form $\zq u$ are sufficiently homogeneous. 
This is not immediately clear here, 
unlike the case of many other homogeneity claims. 
Assume that conditions $p,q\in\Qa$ satisfy 
the next requirement:%
\busm
{*sut*}
{$\pf p=\pf q$ \ \ and \ \ 
$\ps p\cup \ps q\sq \pfv p=\pfv q$.
}
\noi

Then, a transformation $\yh pq$ acting 
on conditions is defined as~follows. 

If $p=q$, then define $\yh pq(r)=r$ for all $r\in\Qa$, 
the~identity. 

Suppose that $p\ne q$. 
Then, $p,q$ are incompatible by \eqref{*sut*} and 
Lemma~\ref{*comp}.
Define  
\mbox{$\dyh pq=\ens{r\in\Qa}{r\leq p\lor r\leq q}$}, 
\index{zzdpq@$\dyh pq$}%
\index{domain!dpq@$\dyh pq$}%
\imar{dyh pq}%
the \rit{domain} of $\yh pq$. 
Let $r\in\dyh pq$. 
We put $\yh pq(r)=r':=\pa{\ps {r'}}{\pf{r'}}$, where 
$\pf{r'}=\pf{r}$ and
\imar{yh pq}%
\index{zzhpq@$\yh pq$}%
\busr
{*sut*1}
{\ps {r'}=
\left\{
\bay{rcl}
(\ps {r}\bez \ps p)\cup \ps q 
&\text{ in case }& r\leq p\,,\\[1ex]
(\ps {r}\bez \ps q)\cup \ps p 
&\text{ in case }& r\leq q\,.\\[1ex]
\eay
\right.}

In this case,  
the difference between $\ps {r}$ and $\ps {r'}$ 
is located within the set $X=\pfv p=\pfv q$, so~that 
$\ps r\cap X=\ps p$ and $\ps {r'}\cap X=\ps q$ 
whenever $r\leq p$, while 
$\ps r\cap X=\ps q$ and $\ps {r'}\cap X=\ps p$ 
whenever $r\leq q$. 
The next lemma is Lemma 6 in \cite{kl57}.

\ble
\label{homL}
\ben
\renu
\jtlb{homL1}\msur%
If\/ $u\sq\Fun$ is dense and\/ 
$p_0,q_0\in\zq u$, 
then there exist conditions\/ $p,q\in\zq u$ with\/ 
$p\leq p_0$, $q\leq q_0$, satisfying \eqref{*sut*}.

\itlb{homL2}%
Let\/ $p,q\in\Qa$ satisfy\/ \eqref{*sut*}. 
If\/ $p=q$, then\/ $\yh pq$ is the identity 
transformation.
If\/ $p\ne q$, then\/ $\yh pq$ is an order automorphism
of\/ $\dyh pq=\ens{r\in\Qa}{r\leq p\lor r\leq q}$, 
satisfying\/ $\yh pq(p)=q$ and\/ 
$\yh pq=(\yh pq)\obr=\yh qp$.  

\itlb{homL3}%
If\/ $u\sq\Fun$ and\/ 
$p\yi q\in \zq u$ satisfy\/ \eqref{*sut*},
then\/ $\yh pq$ maps the set\/ 
$\zq u\cap \dyh pq$ 
onto\/ itself 
order-preserving. 
\een 
\ele

%\begin{proof}[Proof (sketch).]
%
\bpf[sketch]
\ref{homL1}
By the density of $u$, there is a countable set $F\sq\Fun$ 
satisfying 
$\pf p\cup\pf q\sq F$ and   
\mbox{$\ps p\cup \ps q\sq F^{\vee}=
\ens{f\res\xi}{f\in F\land 1\le\xi<\bom}$}. 
Put $p=\pa{\ps p}{F}$ and $q=\pa{\ps q}{F}$.
Claims \ref{homL2} and \ref{homL3}
are routine.
\epf
%\end{proof}

\bcor
[in $\rL$]
\lam{homC}
If a set\/ $u\sq\Fun$ is dense,  
then\/ $\zq u$ is\/ {\em\ubf cone homogeneous} 
in the sense \mbox{of\/ {\rm\cite{dob_sdf}}}, \ie, 
if\/ $p_0,q_0\in\zq u$,
then there exist conditions\/ $p,q\in\zq u$ with\/ 
$p\leq p_0$, $q\leq q_0$, such that the cones\/ 
$\zq u_{\leq p}=\ens{p'\in\zq u}{p'\leq p}$ and\/ 
$\zq u_{\leq q}$ are order-isomorphic.
%\hl{\qed} 
%MDPI: suggest remove it.
%%%%REMOVED
\ecor

\parf{Systems and product almost-disjoint forcing}
\las{*s41}

To prove Theorem~\ref{mt}, we make use of a forcing 
notion equal to the countable-support product of a 
collapse forcing $\dC$ and \dd{\ombl}many forcing 
notions of the form $\zq u$, $u\sq\Fun$. 

{\ubf We work in $\rL$.} 
Define $\dC={\pws\om\cap\rL}^{<\om}$, 
the set of all 
finite sequences of subsets of $\om$ in $\rL$, 
an ordinary forcing $\pws\om\cap\rL$ to 
collapse down to $\om$. 

Let $\pwb=\ombl$ 
and $\pwo=\pwb\cup\ans{-1}$, 
the \rit{index set} of the mentioned product.  
\index{zzIskr@$\pwo$}%
\index{zzIskr-@$\pwb$}%
\imar{pwo,pwb}%
Let a \rit{system} be any map $U:\abs U\to\pws\Fun$ 
such that $\abs U\sq\pwb$, each set 
$U(\inu)$ ($\inu\in\abs U$) is dense in $\Fun$, 
and the \rit{components} 
$U(\inu)\sq\Fun\;\,(\inu\in\abs U)$ 
are pairwise disjoint. 
\index{system}%

Given a system $U$ in $\rL$, 
we let $\jq U$ be the finite-support 
product of $\dC$ and the sets $\zq{U(\inu)}$, 
$\inu\in\abs U$. 
\index{zzQ@$\fQa$}%
\imar{fQa}%
That~is, $\jq U$ consists of all maps $p$ defined on 
a finite set $\dom p=\abp p\sq\abs U\cup\ans{-1}$ so that 
$p(\inu)\in\zq{U(\inu)}$ 
%\bez\ans{\usl\pu\pu}$ 
for all $\inu\in\abm p:=\abp p\bez\ans{-1}$, 
\imar{abm p}%
and if $-1\in\abp p$, then 
$\bij p:=p(-1)\in\dC$. 
\index{zzzzbp@$\bij p$}%
\index{zzzzpabs@$\abp p$}%
\index{zzzzpabs-@$\abm p$}%
If $p\in \jq U$, then put 
\imar{qf p nu}%
$\qf p \inu=\pf {p(\inu)}$ and  
\index{zzspn@$\qs p\inu$}%
\index{zzfpn@$\qf p \inu$}%
%\index{zzfpvn@$\qfv p \inu$}%
\imar{qs p nu}%
$\qs p\inu=\ps {p(\inu)}$ 
for all $\inu\in\abm p$, so that   
$p(\inu)=\usl{\ps p(\inu)}{\pf p(\inu)}$. 

We order $\jq U$ component-wise: $p\leq q$ 
($p$ is stronger as a forcing condition) iff 
$\abp q\sq\abp p$, 
$\bij q\sq \bij p$ in case $-1\in\abp q$, 
and~$p(\inu)\leq q(\inu)$ in $\zq{U(\inu)}$ 
for all $\inu\in\abm q$.  
Note that $\jq U$ contains 
\rit{the empty condition} $\bo\in\jq U$ 
satisfying $\abp\bo=\pu$;  
\index{zzzzodot@$\bo$}%
\index{condition!empty@empty condition $\bo$}%
obviously, $\bo$ is the \dd\leq least 
(and weakest as a forcing condition) 
element of $\jq U$.

\ble
[in $\rL$]
\lam{*ccc2}
If\/ $U$ is a system, then 
the forcing notion\/ $\jq U$ 
satisfies\/ \dd\ombl CC.  
\ele

\bpf
We argue in $\rL$. 
Assume towards the contrary that 
$X\sq\jq U$ is an antichain of cardinality 
$\card X=\omb$. 
As $\card\dC=\omi$, we can assume that $\bij p=\bij q$ 
for all $p,q\in X$. 
Consider the set $S=\ens{\abm p}{p\in X}$; it consists 
of finite subsets of $\omb$. 

{\it Case 1\/}: $\card S\le\omi$. 
Then, by the cardinality argument, there is a set $X'\sq X$ 
and some $a\in S$ such that $\abm p=a$ for all $p\in X'$ 
and still $\card {X'}=\omb$. 
Note that if $p\ne q$ belongs to $X'$, then $\bij p=\bij q$ 
by the above; therefore, as $p,q$ are incompatible, we 
have $\ps p\ne \ps q$. 
Thus, $P=\ens{\ps p}{p\in X'}$ still satisfies $\card P=\omb$. 
This is a contradiction since obviously 
the set $\ens{\ps p}{p\in \jq U\land \abm p=a}$ has 
cardinality $\omi$.

{\it Case 2\/}: $\card S=\omb$. 
Then, by the \dd\Delta system lemma 
(see \eg\ %please add comma following "e.g.," 
Lemma 111.2.6 in Kunen~\cite{kun}) 
there is a set $S'\sq S$ 
and a finite set $\da\sq\omb$ (the root) 
such that $a\cap b=\da$ for all $a\ne b$ in $S'$, 
and still $\card {S'}=\omb$. 
For any $a\in S$, pick a condition $p_a\in X'$ with 
$\abm p=a$; then, $X''=\ens{p_a}{a\in S'}$ still 
satisfies $\card {X''}=\omb$. 
By construction, if $p\ne q$ belong to $X''$, then 
$\abm p\cap\abm q=\da$ and $p,q$ are incompatible; 
hence, the restricted conditions $p\res\da$, $q\res\da$ 
are incompatible as well. 
Thus, the set $Y=\ens{p\res\da}{p\in X''}$ still has 
cardinality $\card Y=\omb$ and is an antichain. 
On the other hand, $\abm q=\da$ for all $q\in Y$. 
Thus we have a contradiction as in Case 1.
\epf

\parf{Jensen--Solovay construction}
\las{*71} 

Our plan is to define a system $\xU\in\rL$ such that 
any \dd{\jq \xU}generic extension of $\rL$ has a 
subextension that witnesses Theorem~\ref{mt}. 
Such a system will be defined in the form of a 
component-wise union of a \dd\ombl long increasing 
sequence of \rit{small} systems, where the smallness 
means that, in $\rL$, 
the system involves only \dd\omil many 
functions in $\Fun$.

{\ubf We work in $\rL$.} 
\bit
\item  
A system $U$ is 
\index{system!small}%
\rit{small}, if~both the set $\abs U$ and each  
set $U(\inu)\;\,(\inu\in\abs U)$ 
has cardinality $\le\omil$. 

\vitem  
If $U,V$ are systems, $\abs U\sq\abs V$,
and $U(\inu)\sq V(\inu)$ for all $\inu\in\abs U$,
%(or simply $U(i)\sq V(i)$ for all $i\in\pwo$, 
%with reference to the identification)
then say that $V$ \rit{extends} $U$, 
in~symbol $U\pce V$.%
\imar{U pce V}%
\index{extends}%
\index{zzzzlecur@$\pce$}%  

\vitem
If $\sis{U_\xi}{\xi<\la}$ is a \dd\pce increasing 
sequence of systems, 
then define a system 
$U=\lis_{\xi<\la}U_\xi$ 
\imar{lis}
by $\abs U=\bigcup_{\xi<\la}\abs{U_\xi}$ and 
$U(\inu)=\bigcup_{\xi<\la,\:\inu\in \abs{U_\xi}}U_\xi(\inu)$ 
for all $\inu\in\abs U$.
\index{zzzzbigveeUxi@$\lis_{\xi<\la}U_\xi$}% 
\eit
 
We let $\zfcm$ be $\zfc$ minus the Power Set axiom, 
\index{zzzzzfc-@$\zfcm$}%
\index{theory!zfc-@$\zfcm$}%
with the 
schema of Collection instead of Replacement,  
with \AC\ in the form of the well-orderability 
principle, and~with the axiom: \lap{$\omi$ exists}.
See~\cite{gitPWS} on versions of $\zfc$ sans 
the Power Set axiom in~detail. 
Let $\zfcd$ be $\zfcm$ plus the axioms: 
$\rV=\rL$, and~the axiom   
\lap{every set $x$ satisfies $\card x\le\omi$}.%
\index{theory!zfc1@$\zfcd$}%

Let $U,V$ be systems. 
Suppose that $M$ is any transitive model of $\zfcd$ 
containing $\bom$. 
\mbox{Define $\pcm U{U'}M$} iff $U\pce U'$ 
\imar{pcm U U' M}%
and the following holds:
\index{zzUMV@$\pcm U{U'}M$}
\ben
\aenu
\jtlb{pcm1}%
the set 
$\raz U{U'}=
\bigcup_{{\inu\in\abs U}} (U'(\inu) \bez U(\inu))$ 
is \rit{multiply\/ $\Seq$-generic} over $M$,  
\index{multiply $\Seq$-generic}%
in the sense that every sequence 
$\ang{f_1,\,\dots f_m}$ of pairwise 
different functions 
$f_\ell\in \raz U{U'}$ 
%(i.e., $k\ne\ell\imp f_k\nT f_\ell$) 
is generic over $M$ 
in the sense of $\Seq={\omi}^{<\omi}$ as the forcing 
notion in $\rL$, \ and

\itlb{pcm2}%
if $\inu\in\abs U$, then the set 
$U'(\inu)\bez U(\inu)$ is dense in $\Fun$, and therefore  
uncountable.\vom\vim
\een
Note a corollary of \ref{pcm1}: $\raz U{U'}\cap M=\pu$.

\bit
\item 
Let $\jsp$, \rit{Jensen--Solovay pairs}, 
\index{zzSJSP@$\cjsp$}%
be the set of all pairs $\pa MU$, where  
$M\mo\zfcd$ is a transitive 
model containing $\bom$ and $U\in M$ 
is a system.  
Then, the sets $\Seq\yd \jq U$ also belong to $M$.\vom\vim

\item
Let $\sjsp$, \rit{small Jensen--Solovay pairs}, 
be the set of all pairs $\pa MU\in\jsp$   
\index{zzJSP@$\jsp$}%
such that $U$ is a small system in the sense above 
and $\card M\le\omi$ 
(in $\rL$).\vom

\item
$\pa MU\pce\pa{M'}{U'}$ 
($\pa{M'}{U'}$ extends $\pa{M}{U}$) \, iff \,
\index{zzzzlecur@$\pce$}%
$M\sq M'$ and 
$\pcm{U}{U'}{M}$;\vom

%\item
\hspace*{-5ex}
$\pa MU\prec\pa{M'}{U'}$ 
(strict) 
\, iff \,
\index{zzzzlcur@$\prec$}%
$\pa MU\pce\pa{M'}{U'}$ and 
$\kaz \inu\in\abm U \,(U(\inu)\sneq U'(\inu))$.\vom

\item
A~\rit{Jensen--Solovay sequence} 
of length $\la\le\bomp=\omb$ is 
\index{Jensen--Solovay sequence}
any strictly \dd\prec increasing \mbox{\dd\la sequence} 
$\sis{\pa{M_\xi}{U_\xi}}{\xi<\la}$ of 
pairs $\pa{M_\xi}{U_\xi}\in\cjsp$, 
satisfying $U_\eta=\lis_{\xi<\eta}U_\xi$ on limit steps. 
%, and 3) $U_0(i)=\pu$ for all $i\in\pwb$. 
Let $\lt\la$ be the set of all such sequences.\vom
\index{zzJSla@$\lt\la$}%

\item
A pair $\pa MU\in\cjsp$ 
\rit{solves} a set 
\index{solves a set}%
$D\sq\cjsp$ iff either $\pa MU\in D$ or there is no pair   
$\pa{M'}{U'}\in D$ that extends $\pa MU$.\vom\vim

\item 
Let $\sol D$ be the set of all pairs $\pa MU\in\cjsp$, 
which 
\index{zzDsolv@$\sol D$}%
solve a given set $D\sq\cjsp$.\vom\vim 

\item
Let $n\ge3$. 
A sequence 
$\sis{\pa{M_\xi}{U_\xi}}{\xi<\omb}\in\lt\omb$ 
is \dd n\rit{complete} iff it
\index{sequence!\dd ncomplete}%  
intersects every set of the form $\sol D$, 
where $D\sq\cjsp$ is a $\hbs{n-2}(\hb)$ set.\vom
\eit
If $\ka$ is a cardinal then 
$\mathrm{H}\ka$ is the collection of all 
sets $x$ whose transitive closure $\TC(x)$ has cardinality 
$\card{(\TC(x))}<\ka$. 
{Arguing in $\rL$, we have $\hb=\rL_{\omb}$, of course.}

%Recall that $\hb$ is the collection of all sets $x$ 
%whose transitive closure $\TC(x)$ has cardinality 
%$\card{(\TC(x))}<\omb$. 
{\sloppy
Further, $\hbs{n-2}(\hb)$ means definability by a 
$\is{}{n-2}$ formula of the \dd\in language, 
in which any definability parameters in $\hb$ are allowed, 
while $\hbs{n-2}$ means the parameter-free definability. 
Similarly, $\hbd{n-1}(\ans{\bom})$ in the next 
theorem means that $\bom=\omil$ 
is allowed as a sole parameter. 
It is a simple exercise that sets 
$\ans\Seq$ and $\Seq$ 
are $\hbd1(\ans\bom)$ under $\rV=\rL$. 
To account for $\bom$ as a parameter, 
note that the set\/ $\omi$ is $\hbs1$; 
hence, the singleton $\ans\omi$ is $\hbd2$. 
}

Generally, we refer to, e.g., 
\cite{skml}, Part B, 5.4,  
or \cite{jechmill}, Chap.~13, 
on the L\'evy hierarchy of $\in$-formulas and 
definability classes $\is H n\yd \ip Hn\yd\id H n$ 
\index{definability classes!%
$\is H n\yd \ip Hn\yd\id H n$}% 
\index{zzSHC@$\is H n,\,\ip Hn,\,\id H n$}% 
for any transitive set $H$.

\bte
%\begin{Theorem}
[Theorem 3 in \cite{kl57}]
\lam{*7euv}
It is true in\/ $\rL$ that if\/ $n\ge2$, then 
there is a sequence\/ 
$\sis{\pa{M_\xi}{U_\xi}}{\xi<\omb}\in\lt\omb$   
of class\/ $\hbd{n-1}(\ans{\bom})$; 
hence, $\hbd{n-1}$ in case\/ $n\ge3$, 
and in addition\/ 
\dd ncomplete in case\/ $n\ge3$,~such that\/ 
$\xi\in\abs{U_{\xi+1}}$ for all\/ $\xi<\omb$.
%\hl{\qed} 
%
\ete

Similar theorems were established in \cite{kl34,kl36,kl38} 
for different purposes.

\bdf
[in $\rL$]
\lam{*fixu}
{\ubf Fix a  number\/ $\nn\ge2$ during the 
proof 
of Theorem~\ref{mt}.}

Let 
$\sis{\pa{\xM_\xi}{\xU_\xi}}{\xi<\omb}\in\lt\omb$
\imar{xM xU}% 
be a Jensen--Solovay sequence as~in 
Thm~\ref{*7euv}, thus  
\index{zzMxi@$\xM_\xi$}%
\index{zzUxi@$\xU_\xi$}%
\ben
\renu
\jtlb{fixu1}
the sequence is of class\/ $\hbd{\nn-1}$; 

\itlb{fixu2}
we have\/ $\xi\in\abs{\xU_{\xi+1}}$ for all\/ $\xi$; 

\itlb{fixu3}
if\/ $\nn\ge3$, then the sequence is\/ \dd\nn complete.\vom
\een
Put $\xU=\lis_{\xi<\omb}\xU_\xi$, 
\imar{xU}
so $\xU(\inu)=
\bigcup_{\xi<\omb,\inu\in\abs{\xU_\xi}}\xU_\xi(\inu)$ 
for all $\inu\in\pwb$. 
Thus, $\xU\in\rL$ is a system and $\abs\xU=\pwb$ 
since $\xi\in\abs{\xU_{\xi+1}}$ for all $\xi$.

We define $\dQ=\jq \xU$ (the {\em\ubf basic forcing notion}). 
\imar{dQ}%
%and $\dQ_\xi=\jq{\xU_\xi}$ for $\xi<\omb$. 
%
\index{zzUb@$\xU$}%
\index{zzPb@$\dQ$}%
\index{zzPxi@$\dQ_\xi$}%
Thus, $\dQ\in\rL$ is the finite-support product of the 
set $\dC$ and sets 
$\dQ(\inu)=\zq{\xU(\inu)}\yt \inu\in\pwb$.
\edf
 
\ble
[in $\rL$]
\lam{*dl}
The binary relation\/ $f\in\xU(\inu)$ 
%and  the set\/ $\dQ$ 
is\/ $\hbd{\nn-1}(\ans{\bom})$. 
\ele

\bpf
To get a $\is{}{}$ definition, 
make use of \ref{fixu1} of Definition~\ref{*fixu}. 
To get a $\ip{}{}$ definition, note that, in $\rL$, 
$f\in\xU(\inu)$ iff for any $\xi<\omb$, if 
$f\in \xM_\xi$ and $\nu<\xi$ then $f\in\xU_\xi(\nu)$.
\epf

\parf{Basic generic extension}
\las{*Pex}

We consider $\dQ_\nn:=\dQ_\nn=\jq\xU$ 
(see Definition~\ref{*fixu}) 
as a forcing notion in $\rL$.
Accordingly, we will study \dd{\dQ}generic 
extensions $\rL[G]$ 
of the ground universe $\rL$.  
Define some elements of these extensions. 
Suppose that $G\sq\dQ$.
%Put $\abm{G}=\bigcup_{p\in{G}}\abm p$; 
%\index{zzGabs@$\abs{G}$}%
\index{zzzzabmG@$\abm G$}% 
%$\abm G\sq\pwb$.
Let
$$
\bij G=\textstyle\bigcup_{p\in G} \bij p\,,
\index{zzbG@$\bij G$}%.
\quad \text{and}\quad
\qs G \inu =\ps {G(\inu)}=
\textstyle\bigcup_{p\in G}\qs p \inu 
$$
\index{zzGi@$G(\inu)$}%
\index{zzSGi@$\qs G \inu$}%
for any $\inu\in\abm G$, where 
$G(\inu)=\ens{p(\inu)}{p\in G}\sq\zq{\xU(\inu)}$.
Thus, $\qs G \inu \sq\Seq$. 

Therefore, 
any \dd{\dQ}generic set $G\sq\dQ$ splits into 
the family of sets $G(\inu)\yt \inu\in\pwb$, 
%and each $G(\inu)$ is \dd{\zp{U(\inu)}}generic, 
and a separate 
%generic 
map $\bij G:\om\onto\pws\om\cap\rL$. 
It follows from Lemma~\ref{*ccc2} by standard arguments  
that \dd{\dQ}generic extensions of $\rL$ 
satisfy $\omi=\om_2^\rL$. 

\ble
[Lemma 9 in \cite{kl57}]
\lam{*krg}
Let\/ 
%$U$ be a system in\/ $\rL$, and\/ 
$G\sq \dQ$ be a set\/ \dd{\dQ}generic 
over\/ $\rL$. 
Then\/$,$\vom\vim
\ben
\renu
\jtlb{krg5}\msur% 
$\bij G$ is a\/ \dd\dC generic map from\/ 
$\om$ onto\/ $\pws\om\cap\rL\;;$\vim

\itlb{Kreg1}%
if\/ $\inu\in\pwb$, then the set\/ 
$G(\inu)=\ens{p(\inu)}{p\in G}\in\rL[G]$ 
is\/ \dd{\zp{\xU(\inu)}}generic over\/ $\rL$ 
--- hence,
if\/ $f\in\Fun$, then\/ 
$f\in \xU(\inu)$ iff\/ 
$\qs G \inu \text{ does not cover }f\;.$
\qed
\een
\ele
Now suppose that $c\sq\pwo$. 
If $p\in\dQ$, then a {\em restriction\/} 
$p'=p\res c\in\dQ$ is defined by 
\index{zzp"1c@$p\res c$}%
$\abm{p'}=c\cap\abm p$ and $p'(\inu)=p(\inu)$ 
for all $\inu\in\abm{p'}$. 
%, and $p'(i)=\usl{\pu}{\pu}$ for $i\nin c$.
In particular, if $\inu\in\pwo$, then let 
$$
p\rn \inu=p\res{(\abp p\bez\ans \inu)}
\qand 
p\rne\inu=p\res{\ans \inu}\;\;
\text{(identified with $p(\nu)$).}
$$
\index{zzp"1ni@$p\rn \inu$}%
If $G\sq\dQ$, then let 
$G\res c=\ens{p\in G}{\abm p\sq c}$
(=$\ens{p\res c}{p\in G}$ 
in case $c\in\rL$).

\index{zzG"1c@$G\res c$}%
Put 
$G\rn \inu=\ens{p\in G}{\inu\nin\abp p}
=G\res{(\pwo\bez\ans \inu)}$.
\index{zzG"1ni@$G\rn \inu$}%
          
Writing $p\res c$, 
%$G\res c$,
it is not assumed that $c\sq\abp p$. 

%\edf

The proof of Theorem~\ref{mt} makes use of a generic 
extension of the form $\rL[G\res c]$, where $G\sq\dQ$ 
is a set 
\dd\dQ generic over $\rL$ and~$c\sq\pwo\yt c\nin\rL$.  
%The following two theorems 
%will play the key role in the proof. 

Define formulas $\dGa_\inu$ ($\inu\in\pwb$) 
as follows: 
%(compare with \S\,4.5 in~\cite{kl56}):
$$
\dGa_\inu(S)\,:=_{\rm def}\;S\sq\Seq\,\land\,
%\kaz j\in\om\;
\kaz f\in \Fun\;
\big(f\in \xU(\inu) \eqv 
{S}\text{ does not cover }f). 
$$ 

\ble
[Lemma 22 in \cite{kl57}] 
\lam{gi}
Suppose that a set\/ 
$G\sq\dQ$ is\/ \dd\dQ generic over\/ $\rL$ and\/ 
$\inu\in\pwb$, $c\in\rL[G]\yt \pu\ne c\sq\pwo$. 
Then, \/ $\omi^{\rL[G\res c]}=\om_2^\rL$ and
\ben
\renu
\jtlb{gi1}\msur%
%if\/ $\inu\in c$ then\/ $S_G(\inu)\in\rL[G\res c]$ and\/ 
$\dGa_\inu(\qs G \inu )$ holds$;$ 
%but 

\itlb{gi2}\msur%
%if\/ $\inu\nin c$ then\/ 
$\qs G \inu \nin\rL[G\rn\inu]$---generally, 
there are no sets\/ $S\sq\Seq$ in\/ 
$\rL[G\rn\inu]$ satisfying\/ 
$\dGa_\inu(S)\,;$ 

\itlb{gi0+}%
if\/ $-1\in c$, then\/ $\bij G\in\rL[G\res c]$, 
and if\/ $\inu\in c$, then\/ $\qs G \inu \in\rL[G\res c]$.
\qed
\een
\ele

The next key theorem is Theorem 4 in \cite{kl57}. 
Note that if $\nn=2$, then the result is an easy 
corollary of the Shoenfield absoluteness theorem. 

\bte
[elementary equivalence theorem]
\lam{*eet}
Assume that in\/ $\rL$, $-1\in d\sq\pwo$, 
sets\/ $Z',Z\sq\pwb\bez d$ satisfy\/ 
$\card{(\pwb\bez Z)} \le \omi$ and\/   
$\card {(\pwb\bez Z')} \le \omi$, 
the symmetric difference\/ $Z\sd Z'$ is at most countable 
and the complementary set\/ 
$\pwb\bez(d\cup Z\cup Z')$ is infinite. 

Let\/ $G\sq\dQ$ be\/ \dd\dQ generic over\/ $\rL$, 
and\/ $x_0\in\rL[G\res d]$ be any real. 

Then,\/ any closed\/ $\is1\nn$ formula\/ $\vpi$, 
with real parameters in\/ $\rL[x_0]$, 
is simultaneously true in the models\/ 
$\rL[x_0,G\res Z]$ and\/ $\rL[x_0,G\res Z']$.
\qed  
\ete

\parf{The model}
\las{175}

Here, we introduce a submodel of the basic \dd\dQ generic 
extension $\rL[G]$ defined in Section~\ref{*Pex} that will 
lead to {the proof}  of Theorem~\ref{mt}.

Recall that a number $\nn\ge2$ is fixed by 
Definition~\ref{*fixu}. 

Under the assumptions and 
notation of Definition~\ref{*fixu},
consider a set $G\sq\dQ$, \dd\dQ generic over $\rL$. 
\mbox{Then, $\bij G=\bigcup G(-1)$} is a 
\index{zzbG@$\bij G$}%
\dd\dC generic map from $\om$ onto $\pws\om\cap\rL$ 
by Lemma~\ref{*krg} \ref{krg5}.
We define%
\index{zzwG@$\wg G$}%
\index{zzWG@$\Wg G$}%
\imar{fn sli}%
\busr{W=}%
{\wg G=
\ens{\om k+2^j}{k<\om\land j\in \bij G(k)}
\cup 
\ens{\om k+3^j}{j,k<\om}\sq\om^2,}
%\eus%
%
and $\Wg G=\ans{-1}\cup \wg G$.
We also define, for~any $m<\om$,  
$$
\index{zzwmG@$\wgb m G$}%
\index{zzwmG@$\wgm m G$}%
\wgb mG=\ens{\om k+\ell\in \wg G}{k\ge m}\,,
\quad 
\wgm mG=\ens{\om k+\ell\in \wg G}{k<m}\,,
$$
and accordingly, 
\index{zzWmG@$\Wgb m G$}%
\index{zzWmG@$\Wgm m G$}%
$\Wgb mG=\ans{-1}\cup \wgb mG$ and 
$\Wgm mG=\ans{-1}\cup \wgm mG$.

With these definitions, each $k$th slice 
\busr
{wk=}
{\wgk k G=
\ens{\om k+2^j}{j\in \bij G(k)}
\cup 
\ens{\om k+3^j}{j<\om}
\index{zzwkG@$\wgk k G$}%
}
of $\wg G$ is 
necessarily infinite and coinfinite, and~it 
codes the target 
set $\bij G(k)$ since 
\busr
{bjk}
{\bij G(k)=\ens{j<\om}{\om k+2^j\in \wgk k G}
=\ens{j<\om}{\om k+2^j\in \Wg G}.}

Note that definition 
\eqref{W=} is 
\rit{monotone \poo\ $\bij G$}, \ie, 
if $\bij G(k)\sq \bij{G'}(k)$ for all $k$, 
then $\wg G\sq\wg {G'}$ and 
$\Wg G\sq\Wg {G'}$.
Anyway, $\wg G\sq\om^2$ (the ordinal product) 
is a set in the model 
$\rL[\bij G]=\rL[\Wg G]=\rL[\wg G]=\rL[\wgb mG]$ 
for each $m$, whereas $\wgm mG\in\rL$ 
for all $m$. 
Finally, let 
$ 
W=[\om^2,\omb)=\ens{\za}{\om^2\le\za<\omb}\,.
\index{zzW@$W=[\om^2,\omb)$}%
$ 

Recall that if $c\sq\pwo$, then 
$G\res c=\ens{p\in G}{\abp p\sq c}$. 

\bce
\ubf
To prove Theorem~\ref{mt}, 
we consider the model \ 
$\rL[G\res (\Wg G\cup W)]\,\sq\,\rL[G]$.
\ece

\bte
\lam{33}
If\/ $G$ is 
a\/ \dd\dQ generic set over\/ $\rL$, then 
the class\/ $\rL[G\res (\Wg G\cup W)]$ 
suffices to prove Theorem \ref{mt}. 
That is, \lsep\ holds 
in\/ $\rL[G\res (\Wg G\cup W)]$ both for\/   
$\is1{\nn+1}$ sets of integers and for\/   
$\ip1{\nn+1}$ sets of integers. 
\ete 

The proof will include several lemmas. 

For the next lemma, 
we let $\for_\dQ$ be the $\dQ$-forcing notion 
defined in $\rL$. 
If $p\in\dQ$ and $-1\in\abp p$, 
then let $p\rne{-1}:=p\res\ans{-1}$. 
This can be identified with just $p(-1)\in\dC$, of 
course, but formally $p\rne{-1}\in\dQ$. 
If $-1\nin\abp p$, then let $p\rne{-1}:=\bo$ 
(the empty condition). 
Let $\naG$ be the canonical \dd\dQ name 
for the generic set $G\sq\dQ$, 
$\check W$ be a name for the set
$W=[\om^2,\ombl)\in\rL$, and 
$\bijn$ be a canonical 
\dd\dQ name for $\bij G$.
\index{name!Gund@$\naG$}%
\index{zzzzGund@$\naG$}%

%forced to be true in $\rL[\dgwgw]$ 

\ble
[reduction to the $\dC$-component]
\lam{-1}
Let\/ $p\in\dQ$ and let\/ $\Phi(\bijn)$ 
be a closed formula containing only\/ 
$\bijn$ and names for sets in\/ $\rL$ 
as parameters. 
Assume that\/ 
$$
p\for_\dQ 
``\Phi\text{ is true in }\rL[\dgwgw]''\,.
$$ 

Then,\/ $p\rne{-1}\for_\dQ 
``\Phi\text{ is true in }\rL[\dgwgw]''$ 
as well. 
\ele

\bpf
By the product forcing theorem, if $G\sq\dQ$ is 
\dd\dQ generic over $\rL$, then the model 
$\rL[\gwgw]$ is a 
\dd{\dQ'}generic extension of $\rL[\bij G]$, where 
$\dQ'=
\prod_{\inu\in\Wg{G}\cup W}\zq{\xU(\inu)}$ 
is a forcing in $\rL[\bij G]$. 
However, it follows from Corollary~\ref{homC} that 
$\dQ'$ is a (finite-support) product of 
cone-homogeneous forcing notions. 
Therefore, $\dQ'$ itself is a cone homogeneous forcing, 
and we are finished (see e.g., Lemma 3 in \cite{dob_sdf} 
or Theorem IV.4.15 in \cite{kun}).
\epf

\parf{Two key lemmas}
\las{2key}

{The following two lemmas present two key properties 
of models of the form $\rL[G\res (\Wg G\cup W)]$ 
involved in the proof of Theorem~\ref{33}. 
The first lemma 
%, which we take from \cite{kl57}, 
shows that all constructible reals are 
$\id1{\nn+1}$ in such a model.}

\ble
%[Lemma 23(i) in \cite{kl57}]
\lam{7lem1}
Let a set\/ $G\sq\dQ$ be\/ 
\dd\dQ generic over\/ $\rL$. 
Then, it holds in\/ $\rL[G\res (\Wg G\cup W)]$ that\/ 
$\wg G$ is $\is1{\nn+1}$ and  
each set\/ $x\in\rL\yt x\sq\om$ is\/ $\id1{\nn+1}\;.$ 
\ele

\bpf
Consider an arbitrary ordinal $\inu=\om k+\ell$; 
$k,\ell<\om$. 
We claim that 
\busr
{gx}
{{\inu\in \wg G}\eqv {\sus S\,\dGa_{\inu}(S)}  
}
holds in $\rL[G\res (\Wg G\cup W)]$. 
Indeed, assume that $\inu\in \wg G$. 
Then, $S=\qs G \inu \in \rL[G\res \Wg G]$, and 
we have $\dGa_{\inu}(S)$ in $\rL[G\res (\Wg G\cup W)]$
by Lemma~\ref{gi} \ref{gi0+},\ref{gi1}.
Conversely, assume that $\inu\nin \wg G$. 
%Then accordingly $\inu=\om k+2^j\nin \Wg G$. 
Then,~we have 
$\Wg G\in\rL[\bij G]\sq \rL[G\res \Wg G]\sq\rL[G\rn\inu]$, 
but $\rL[G\rn\inu]$ contains no $S$ with 
$\dGa_{\inu}(S)$ by Lemma~\ref{gi}~\ref{gi2}.

\nek{\hk} {{\mathrm{H}\bbkappa}}
\nek{\hi} {{\mathrm{H}\omi}}

Let 
$\bbkappa:=\omb^\rL$ 
in the remainder of the proof of the lemma. 
Then it follows that 
$\bbkappa=\omi^{\rL[G\res (\Wg G\cup W)]}=\omi^{\rL[G]}$, 
by {Lemma~\ref{gi}.}

Note that the right-hand side of \eqref{gx} defines a 
$\is{\hk}{\nn}(\ans{\oli,\Seq})$  relation in 
the model $\rL[G\res (\Wg G\cup W)]$ by Lemma~\ref{*dl}.
(Indeed, we have 
%$\omi^{\rL[G\res W_G]}=\omi^{\rL[G]}=\omb^\rL$ and 
$(\hk)^\rL=(\hb)^\rL=\rL_{\omb^\rL}=\rL_{\bbkappa}$
%=(\rL)^\hc$ 
in $\rL[G\res (\Wg G\cup W)]$, therefore 
$(\hk)^\rL$ is $\is\hk1$ in $\rL[G\res(\Wg G\cup W)]$.) 
On the other hand, the~sets $\ans{\oli}$ and 
$\ans{\Seq}$ 
remain $\hbd2$ singletons in $\rL[G\res(\Wg G\cup W)]$; 
they can be eliminated since $\nn\ge2$. 
This yields $\wg G\in\is\hk\nn=\is\hi \nn$ in 
$\rL[G\res(\Wg G\cup W)]$.
It follows that $\wg G\in\is1{\nn+1}$ in 
$\rL[G\res(\Wg G\cup W)]$ by Lemma 25.25 in 
\cite{jechmill}, as~required.

Now, let $x\in\rL\yt x\sq\om$.
By genericity, there exists $k<\om$ such that 
$\bij G(k)=x$. 
Then, $x=\ens{j}{\om k+2^j\in \wg G}$ by \eqref{W=}; 
therefore, $x$ is $\is1{\nn+1}$ as well.
However, $\om\bez x\in\is1{\nn+1}$ by the same argument. 
Thus, $x$ is $\id1{\nn+1}$ in $\rL[G\res(\Wg G\cup W)]$, 
as required. 
\epf

The proof of the next lemma 
%, also from \cite{kl57}, 
%is of the same category as Theorem~\ref{*eet} above, 
%and its proof 
involves Theorem~\ref{*eet} above as a key reference. 
The lemma holds for $\nn=2$ by Shoenfield.

\ble
%[Lemma~24 in~\cite{kl57}]
\lam{sn+7}
Suppose that\/ $G\sq\dQ$ is\/ 
\dd\dQ generic over\/ $\rL$, 
$m<\om$, $c\sq \wgm m G$, $c\in\rL$. 
Then, any closed\/ $\is1\nn$ formula\/ $\Phi$, with~reals in\/ $\rL[G\res (c\cup \Wgb mG\cup W)]$ 
as parameters, 
is simultaneously true 
in\/ $\rL[G\res {(c\cup\Wgb mG\cup W)}]$ 
and in\/ $\rL[G\res {(\Wg G\cup W)}]$.

It follows that if\/ $c'\sq c\sq \wgm m G$ in\/ $\rL$, 
then any closed\/ $\is1{\nn+1}$ formula\/ $\Psi$, with~parameters in\/ $\rL[G\res (c'\cup \Wgb mG\cup W)]$, 
true in\/ $\rL[G\res {(c'\cup\Wgb mG\cup W)}]$, 
is true in the model\/ $\rL[G\res {(c\cup\Wgb mG\cup W)}]$ 
as well.
\ele

\bpf
There is an ordinal $\xi<\omb$ such that 
all parameters in $\vpi$ belong to 
$\rL[G\res Y]$, \mbox{where 
$Y=c\cup \Wgb mG\cup X$} and 
$X=[\om^2,\xi)=\ens{\ga}{\om^2\le\ga<\xi}$. 
%, and accordingly $y$ belongs to 
%$\rL[G\res {(\Wg G\cup X)}]$.
%
The set $Y$ 
belongs to $\rL[\bij G]$; 
in fact, $\rL[Y]=\rL[\bij G]$. 
Therefore, $G\res Y$ 
is equi-constructible with the pair 
$\ang{\bij G,\sis{\qs G \inu }{\nu\in Y}}$. 
Here, $\bij G$ is 
a map from $\om$ onto $\pws\om\cap\rL$. 
It follows that there is a real $x_0$ with 
$\rL[G\res Y]=\rL[x_0]$.  
Then, all parameters of $\vpi$ belong to $\rL[x_0]$. 
%while $y$ belongs to $\rL[x_0,G\res e]$, 
%where $e=\wgm m G\bez c$.

To prepare an application of Theorem~\ref{*eet} 
of Section~\ref{*Pex}, we put 
$$
\bay{rclcccc}
Z' &=& [\xi,\omb)\;,\\[0.5ex]  
Z &=& e\cup Z'\;,\quad\text{where}\quad 
e=\wgm m G\bez c\;,\\[0.5ex]  
d &=& \ans{-1}\cup 
\ens{\om k+j}{k\ge m\land j<\om}\cup X\,.\vim
\eay
$$

It is easy to check that all requirements of 
Theorem~\ref{*eet} for these sets are fulfilled.
Moreover, as 
$\Wgb mG\sq\ans{-1}\cup\ens{\om k+j}{k\ge m\land j<\om}$, 
we have $Y=c\cup \Wgb mG\cup X\sq d$;~hence, $x_0\in\rL[G\res d]$.  
%$y\in\rL[x_0,G\res e]\sq\rL[x_0,G\res Z]$. 
Therefore, we conclude by Theorem~\ref{*eet} 
that the formula  
$\vpi$ is simultaneously true 
in $\rL[x_0,G\res Z]$ and in $\rL[x_0,G\res Z']$. 
However, 
$$
\rL[x_0,G\res Z']=\rL[G\res (Y\cup Z')]=
\rL[G\res {(c\cup \Wgb mG\cup W)}]
$$ 
by construction, 
while~\mbox{$\rL[x_0,G\res Z]=\rL[G\res {(\Wg G\cup W)}]$}, 
and we are done.
\epf

\parf{Finalization: \boldmath{$\is1{\nn+1}$} case}
\las{175*} 

Here, we finalize the proof of Theorems~\ref{33} 
and \ref{mt} 
\poo\ $\is1{\nn+1}$ sets of integers. 
We generally follow the line of arguments 
sketched by Harrington ([Addendum A3]~\cite{h74}) 
for the $\is13$ case, 
with suitable changes mutatis mutandis. 
We will fill in all details omitted in \cite{h74}.

Recall that a number $\nn\ge2$ is fixed by 
Definition \ref{*fixu}. 
We assume that\vom 
\ben
\fenu
\jtlb{0xy}% 
a set\/ $G\sq\dQ$ is \dd\dQ generic over $\rL$, 
sets $x,y\sq\om$ belong to $\rL[G\res (\Wg G\cup W)]$,  
and it holds in\/ $\rL[G\res (\Wg G\cup W)]$ that  
$x,y$ are disjoint\/ $\is1{\nn+1}$ sets.\vom
\een

The goal is to prove that $x,y$ can be separated by 
a set $Z\in\rL$ and then argue that $Z$ is $\id1{\nn+1}$ 
by Lemma \ref{7lem1}.
Recall that 
$W=[\om^2,\ombl)=\ens{\xi}{\om^2\le\xi<\ombl}$.
Suppose that\vom
\ben
\fenu
\atc
\jtlb{1stt}\msur
$\vpi(\cdot)$ and $\psi(\cdot)$ are parameter-free 
$\is1{\nn+1}$ formulas that provide   
$\is1{\nn+1}$ definitions for the sets, respectively,\ $x,y$ 
of \ref{0xy}---that is, 
$$
x=\ens{\ell<\om}{\vpi(\ell)}
\qand
y=\ens{\ell<\om}{\psi(\ell)}
$$ 
in $\rL[G\res(\Wg G\cup W)]$. 
The assumed implication 
$\kaz \ell\:(\vpi(\ell)\imp\neg\:\psi(\ell))$ 
(as $x\cap y=\pu$) is 
forced to be true in $\rL[\dgwgw]$ 
by a condition $p_0\in G$.\vom 
\een

Here, $\naG$ is the canonical \dd\dQ name 
for the generic set $G\sq\dQ$ 
while $\check W$ is a name for the set
$W=[\om^2,\ombl)\in\rL$.
\index{name!Gund@$\naG$}%
\index{zzzzGund@$\naG$}%

We observe that 
$\kaz \ell\:(\vpi(\ell)\imp\neg\:\psi(\ell))$ 
is a parameter-free sentence. 
Therefore, it can \noo\ be assumed that 
$\abp{p_0}=\ans{-1}$, by Lemma~\ref{-1}. 
In this case, the condition 
$p_0\in\dQ$ can be identified with its only 
nontrivial component $s_0=p_0(-1)\in\dC$. 

\ble
[interpolation lemma]
\lam{1122}
Under the assumptions of\/ \ref{1stt},
if\/ $\ell,\ell'<\om$, conditions\/ $p,p'\in\dQ$   
satisfy\/ $p\le p_0$ and\/ $p'\le p_0$, 
and we have 
\bce
$p\for_\dQ  ``\rL[\dgwgw]\mo\vpi(\ell)$'' 
\ \ and\/ \ \ 
$p'\for_\dQ  ``\rL[\dgwgw]\mo\psi(\ell')$''. 
\ece
Then,\/ $\ell\ne\ell'$. 
\ele

%\begin{proof}
%[Proof ({\rm sketched in ([Addendum A3]~\cite{h74}) for $\nn=2$}) ]

\bpf
[sketched in {\cite[Addendum A3]{h74}} for $\nn=2$]
First of all, by Lemma~\ref{-1}, we can \noo\ 
assume that $\abp{p}=\abp{p'}=\ans{-1}$; 
so, the components $s=p(-1)$ and $s'=p'(-1)$ satisfy 
$s_0\sq s$ and $s_0\sq s'$ in $\dC$. 

We \noo\ assume that the tuples $s,s'$ have the 
same length $\lh s=\lh{s'}=m$. 
(Otherwise, extend the shorted one by a sufficient 
number of new terms equal to $\pu$.) 
Define another condition $t\in\dC$ such that 
$\dom t=m$ and $t(j)=s(j)\cup s'(j)$ for all $j<m$.
Accordingly, define $q\in\dQ$ so that $\abp q=\ans{-1}$ 
and $q(-1)=t$. 
Despite that $q$ may well be incomparable with $p,p'$ 
in $\dQ$, we claim that 
\bul{ql}
q
\for_\dQ 
\text{\rm``$\rL[\dgwgw]\mo\vpi(\ell)\land \psi(\ell')$''}\,. 
\eus

To prove the \dd\vpi part of \eqref{ql}, let $H\sq\dQ$ 
be a set \dd\dQ generic over $\rL$, and $q\in H$. 
Then, $t\su\bij H$. 
We have to prove that $\vpi(\ell)$ holds in $\rL[\hwhw]$.

Define another generic set $K\sq\dQ$, 
slightly different from $H$, 
so that
\ben
\Aenu
\jtlb{kh1}%
$K(\inu)=H(\inu)$ for all $\inu\in\pwb=\ombl$;\vim
 
\itlb{kh2}%
$s\su\bij K$; \ \ and\vim 

\itlb{kh3}%
 if $m\le j<\om$, then $\bij K(j)=\bij H(j)$.
\een
\noi 

In other words, the only difference between $K$ and $H$ 
is that $\bij K\res m=s$ but $\bij H\res m=t$. 

It follows that $p\in K$; hence, 
$\vpi(\ell)$ holds in $\rL[\kwkw]$ by the assumptions of 
the lemma. 
Now, we note that by definition,
$$
\Wg{K}\cup W=\wgm m{K}\cup \Wgb m{K}\cup W,
\quad 
\Wg{H}\cup W=\wgm m{H}\cup \Wgb m{H}\cup W,
$$

Here, the sets $c_H=\wgm m{H}$ and $c_K=\wgm m{K}$ 
satisfy $c_K\sq c_H$ 
(since $\bij K(j)=s(j)\sq t(j)=\bij H(j)$ for all 
$j<m$). 
In addition, $\Wgb m{H}=\Wgb m{K}$ 
(since $\bij K(j)\bij H(j)$ for all $j\ge m$). 
To conclude, 
\bul{kh}
\Wg{K}\cup W=c_K\cup \Wgb m{H}\cup W,
\quad 
\Wg{H}\cup W=c_H\cup \Wgb m{H}\cup W,
\eus
and $c_K\sq c_H=\wgm m{H}$. 
On the other hand, 
it follows %\hl 
{from} 
(A) that 
$K\res c=H\res c$ for any $c\sq\pwb$, whereas $\bij K$ 
and $\bij H$ are recursively reducible to each other 
by (B),(C). 
Therefore,
$$
\rL[\kwkw]
=\rL[H\res(\Wg{K}\cup W)]
=\rL[H\res(c_K\cup \Wgb m{H}\cup W)]
$$
by \eqref{kh}.
However, $\vpi(\ell)$ holds in this model by the above. 
It follows by Lemma~\ref{sn+7} that $\vpi(\ell)$ 
holds in $\rL[\hwhw]=\rL[c_H\cup \Wgb m{H}\cup W]$ 
as well. 
(Harrington circumvents \mbox{Lemma~\ref{sn+7} in \cite{h74}} 
by a reference to the Shoenfield absoluteness theorem.)
We are finished. 

After \eqref{ql} has been established, we recall that 
$q\le p_0$ in $\dQ$ by construction. 
It follows that $\ell\ne\ell'$ by the choice of $p_0$ 
(see \ref{1stt} above).
\epf

\bpf[Theorems~\ref{33} 
and~\ref{mt}:~$\is1{\nn+1}$ case]
%\bpf[Proof of Theorems~\ref{33} and~\ref{mt}: 
%$\is1{\nn+1}$ case]
We work under the assumptions of \ref{0xy} and \ref{1stt} 
above. 
Consider the following sets in $\rL$:
$$
\bay{rcl}
Z_x
&=&\ens{\ell<\om}
{\sus p\in\dQ\,(p\le p_0\land 
p\for_\dQ  \text{``$\rL[\dgwgw]\mo\vpi(\ell)$''})};\\[1ex]
Z_y
&=&\ens{\ell'<\om}
{\sus p'\in\dQ\,(p'\le p_0\land 
p'\for_\dQ  \text{``$\rL[\dgwgw]\mo\psi(\ell')$''})}.
\eay
$$

Note that $Z_x\cap Z_y=\pu$ by Lemma~\ref{1122}. 
On the other hand, it is clear that $x\sq Z_x$ and 
$y\sq Z_y$ by \ref{1stt}. 
Thus, either of the two sets $Z_x, Z_y\in\rL$ separates 
$x$ from $y$. 
It remains to apply Lemma~\ref{7lem1}. 
\epf

\parf{Finalization: \boldmath{$\ip1{\nn+1}$} case}
\las{176} 

%Here we carry out the proof of Theorems~\ref{33} 
%and~\ref{mt} \poo\ $\ip1{\nn+1}$ sets of integers. 
This will be a mild variation of the argument 
presented in the previous section. 

%\begin{proof}[Proof of Theorems~\ref{33} and~\ref{mt}: 
%$\ip1{\nn+1}$ case, sketch]

\bpf
[Theorems~\ref{33} and~\ref{mt}: $\ip1{\nn+1}$ case, sketch]
Emulating \ref{0xy} and \ref{1stt} above, we 
assume that a set $G\sq\dQ$ is \dd\dQ generic 
over $\rL$, and  
$x,y\sq\om$ are disjoint $\ip1{\nn+1}$ sets in 
$\rL[G\res (\Wg G\cup W)]$, defined by  
parameter-free 
$\ip1{\nn+1}$ formulas, respectively,\ 
$\vpi(\cdot)$ and $\psi(\cdot)$. 
The implication 
$\kaz \ell\:(\vpi(\ell)\imp\neg\:\psi(\ell))$ is 
forced to hold in $\rL[\dgwgw]$ 
by a condition $p_0\in G$ 
satisfying $\abp {p_0}=\ans{-1}$. 
The proof of Lemma~\ref{1122} goes on for 
$\ip1{\nn+1}$ formulas $\vpi,\psi$
the same way, with the only difference that we 
define $t(j)=s(j)\cap s'(j)$ for $j<m$. 
Yet, this is compatible with the application of 
Lemma~\ref{sn+7} because now, $\vpi,\psi$ are 
$\ip1{\nn+1}$ formulas.
\epf

%%%%%%%%%%%%%%%%%%%%%%%%%%%%%%%%%%%%%%%%%%
\parf{Conclusions and discussion}
\las{cd}

In this study, the method of almost-disjoint forcing 
was employed to the problem of obtaining
a model of $\ZFC$ in which the \sep\ principle holds 
for lightface classes $\ip1{n+1}$ and $\is1{n+1}$, 
for a given $n\ge2$, 
for sets of integers.
The problem of obtaining such models has been generally 
known since the early years of modern set theory, 
see, \eg, problems
3029 and 3030 in a survey \cite{mathiasS} by Mathias. 
Harrington \cite[Addendum A3]{h74} claimed the existence of 
such models; yet, a detailed proof has never appeared.

{From our study, it is concluded that the 
%hidden invariance technique (as outlined in Section 6.1)
technique developed in our earlier paper \cite{kl57} 
solves the general case of the problem 
(an arbitrary $n\ge2$) 
by providing a generic extension of $\rL$ in which 
the \lsep\ principle holds 
for classes $\ip1{n+1}$ and $\is1{n+1}$, 
for a given $n\ge2$, 
for sets of integers, for a chosen value $n\ge2$. }

From this result, 
we immediately come to the following problem:

%\bque
%\lam{prob1}
\begin{vopr}
\lam{prob1}
Define a generic extension of\/ $\rL$ in which 
the\/ \lsep\ principle holds 
for classes\/ $\ip1{n+1}$ and\/ $\is1{n+1}$, 
for\/ {\ubf all} $n\ge2$, for sets of integers. 
%\eque
\end{vopr}

The intended solution is expected to be obtained 
on the basis of a suitable product of the forcing notions 
$\dQ_\nn$, $\nn\ge2$, defined in Section~\ref{*Pex}.

And we recall the following major problem.

%\bque

\begin{vopr}
\lam{prob2}
Given\/ $n\ge2$, define a generic extension of $\rL$ 
in which the\/ \sep\ principle holds 
for the classes\/ $\is1{n+1}$ and\/ $\fs1{n+1}$
for\/ {\ubf sets of reals}. 
%\eque
\end{vopr}

The case of \rit{sets of reals} in the \sep\ problem 
is more general, and obviously much 
more difficult, than the case sets of integers.

%\back
%
%The authors are thankful to the anonymous referee
%for a number of helpful and important
%remarks and suggestions, that allowed to substantially 
%improve the text.
%\eack

\bibliographystyle{plain}
\addcontentsline{toc}{subsection}{\hspace*{5.5ex}References}
{\small
\bibliography{65.bib,65kl.bib}%
}

\end{document}